\documentclass[12]{article}

\input epsf
\title{Counting $1$-vertex Triangulations\\ Of Oriented Surfaces}
\author{Roland Bacher\\ Institut Fourier, UMR 5582 CNRS-UJF Universite Grenoble I, \\ BP 74, 38402 StMartin D'Heres, France \\ Roland.Bacher@ujf-grenoble.fr
\and Alina Vdovina \\ Department of Mathematical Sciences,\\ SUNY, Binghamton, NY 13902, USA \\ vdovina@mpim-bonn.mpg.de}
\date{}

\begin{document}
\maketitle

\begin{abstract}
 A {\em $1-$vertex triangulation\/} of an oriented compact
 surface $S$ of genus $g$ is an embedded graph $T\subset S$
 with a unique vertex such that all connected components of 
 $S\setminus T$ are triangles (adjacent to exactly $3$ edges of $T$).

 This paper gives formulas enumerating such triangulations
 (up to equivalence) on an oriented surface of given genus.

 {\em Une triangulation \`a un sommet\/} d'une surface orient\'ee
 compacte $S$ de genre $g$ est un graphe $T\subset S$ qui a un
 unique sommet et dont toutes les faces (composantes connexes de
 $S\setminus T$) sont des triangles (incidentes \`a trois ar\^etes de $T$).

 Cet article donne des formules permettant d'\'enum\'erer
 ces triangulations.
\end{abstract}

\section*{Introduction}

 {\bf Definition 0.1.} 
 A {\em $1-$vertex triangulation\/} of
 an oriented compact surface $S$ of genus $g$ is an embedded
 graph $T\subset S$ with only one vertex such that all connected
 components of $S\setminus T$ are adjacent to exactly $3$ edges
 of $T$ (i.e. are triangles).

 Two such triangulations $T\subset S$ and $T'\subset S'$ are
 {\em isomorphic\/} (or {\em equivalent\/}) if there exists an
 orientation-preserving homeomorphism $\varphi:S \longrightarrow S'$
 such that $\varphi(T)=T'$.

 The aim of this paper is to give formulas for the number of
 such triangulations (up to equivalence) on an oriented surface of
 given genus (such triangulations exist in every genus $g\geq 1$).

 $1-$vertex triangulations have several applications:

 L.Mosher \cite{[M]} has constructed a complex whose fundamental
 group is the mapping class group of an orientable genus $g$ surface. 
 $1-$vertex triangulations appear as the vertices of this complex.

 There is also a bijection between $1-$vertex triangulations
 and Eulerian paths in cubic graphs as considered by 
 J.Brenner and R.Lyndon in \cite{[BL]}. 
 They formulated the problem of classification and enumeration
 of these objects, which is solved in the present paper.
 Brenner and Lindon conjectured, that the automorphism
 group of a $1-$vertex triangulation is always cyclic of order 1,2,3 or
 6. G.Bianchi and R.Cori proved this in a more general form
 in \cite{[BC]}.

 Notice that Brenner and Lyndon considered
 such triangulations from a combinatorial point of view
 motivated by the study of non-parabolic subgroups in the modular group 
 \cite{[BL]}.

 Section \ref{sec:1} introduces oriented Wicks forms (cellular
 decompositions with only one face of oriented surfaces), our main tool.
 Wicks forms are canonical forms for products of commutators
 in free groups \cite{[V]}. Oriented maximal
 Wicks forms and $1-$vertex triangulations are in bijection. 
 Our main theorem will be expressed in the language of Wicks forms.

 Section \ref{sec:2} contains a few facts concerning oriented maximal
 Wicks forms.

 Section \ref{sec:3} contains the proof of our main results.

\section{Main results}
\label{sec:1}

 The objects considered in this section are dual to 1-vertex
 triangulations. They are slightly easier to handle since they carry
 some combinatorial structures more immediately.

 {\bf Definition 1.1.}
 An {\it oriented Wicks form\/} is a cyclic word $w= w_1w_2\dots w_{2l}$
 (a cyclic word is the orbit of a linear word under cyclic permutations)
 in some alphabet $a_1^{\pm 1},a_2^{\pm 1},\dots$ of letters
 $a_1,a_2,\dots$ and their inverses $a_1^{-1},a_2^{-1},\dots$ such that
\begin{itemize}
\item[(i)] if $a_i^\epsilon$ appears in $w$ (for $\epsilon\in\{\pm 1\}$)
 then $a_i^{-\epsilon}$ appears exactly once in $w$,
\item[(ii)] the word $w$ contains no cyclic factor (subword of
 cyclically consecutive letters in $w$) of the form $a_i a_i^{-1}$ or
 $a_i^{-1}a_i$ (no cancellation),
\item[(iii)] if $a_i^\epsilon a_j^\delta$ is a cyclic factor of $w$ then
 $a_j^{-\delta}a_i^{-\epsilon}$ is not a cyclic factor of $w$ 
 (substitutions of the form
 $a_i^\epsilon a_j^\delta\longmapsto x,
 \quad a_j^{-\delta}a_i^{-\epsilon}\longmapsto x^{-1}$ are impossible).
\end{itemize}

 An oriented Wicks form $w=w_1w_2\dots$ in an alphabet $A$ is
 {\em isomorphic\/} to $w'=w'_1w'_2$ in an alphabet $A'$ if
 there exists a bijection $\varphi:A\longrightarrow A'$ with
 $\varphi(a^{-1})=\varphi(a)^{-1}$ such that $w'$ and
 $\varphi(w)=\varphi(w_1)\varphi(w_2)\dots$ define the
 same cyclic word.

 An oriented Wicks form $w$ is an element of the commutator subgroup
 when considered as an element in the free group $G$ generated by
 $a_1,a_2,\dots$. We define the {\em algebraic genus\/} $g_a(w)$ of
 $w$ as the least positive integer $g_a$ such that $w$ is a product
 of $g_a$ commutators in $G$.

 The {\em topological genus\/} $g_t(w)$ of an oriented Wicks
 form $w=w_1\dots w_{2e-1}w_{2e}$ is defined as the topological
 genus of the oriented compact connected surface obtained by
 labeling and orienting the edges of a $2e-$gone (which we
 consider as a subset of the oriented plane) according to
 $w$ and by identifying the edges in the obvious way.

{\bf Proposition 1.1.}
{\sl The algebraic and the topological genus of an oriented Wicks
 form coincide (cf. \cite{[C],[CE]}).}

 We define the {\em genus\/} $g(w)$ of an oriented
 Wicks form $w$ by $g(w)=g_a(w)=g_t(w)$.

 Consider the oriented compact surface $S$ associated to an oriented 
 Wicks form $w=w_1\dots w_{2e}$. This surface carries an immerged graph
 $\Gamma\subset S$ such that $S\setminus \Gamma$ is an open polygon
 with $2e$ sides (and hence connected and simply connected).
 Moreover, conditions (ii) and (iii) on Wicks form imply that $\Gamma$ 
 contains no vertices of degree $1$ or $2$ (or equivalently that the
 dual graph of $\Gamma\subset S$ contains no faces which are $1-$gones
 or $2-$gones). This construction works also
 in the opposite direction: Given a graph $\Gamma\subset S$
 with $e$ edges on an oriented compact connected surface $S$ of genus $g$
 such that $S\setminus \Gamma$ is connected and simply connected, we get
 an oriented Wicks form of genus $g$ and length $2e$ by labeling and 
 orienting the edges of $\Gamma$ and by cutting $S$ open along the graph
 $\Gamma$. The associated oriented Wicks form is defined as the word
 which appears in this way on the boundary of the resulting polygon
 with $2e$ sides. We identify henceforth oriented Wicks
 forms with the associated immerged graphs $\Gamma\subset S$,
 speaking of vertices and edges of oriented Wicks form.

 The formula for the Euler characteristic
 $$\chi(S)=2-2g=v-e+1$$
 (where $v$ denotes the number of vertices and $e$ the number
 of edges in $\Gamma\subset S$) shows that
 an oriented Wicks form of genus $g$ has at least length $4g$
 (the associated graph has then a unique vertex of degree $4g$
 and $2g$ edges) and at most length $6(2g-1)$ (the associated
 graph has then $2(2g-1)$ vertices of degree three and
 $3(2g-1)$ edges).

 We call an oriented Wicks form of genus $g$ {\em maximal\/} if it has
 length $6(2g-1)$. Oriented maximal Wicks forms are dual to 1-vertex
 triangulations. This can be seen by cutting the oriented surface $S$ 
 along $\Gamma$, hence obtaining a polygon $P$ with $2e$ sides. 
 We draw a star $T$ on $P$ which joins an interior point of $P$
 with the midpoints of all its sides. Regluing $P$ we recover $S$
 which carries now a 1-vertex triangulation given by $T$ and each
 1-vertex triangulation is of this form for some oriented maximal
 Wicks form (the immerged graphs $T\subset S$ and $\Gamma\subset S$
 are dual to each other: faces of $T$ correspond to vertices of
 $\Gamma$ and vice-versa. Two faces of $T$ share a common edge if
 and only if the corresponding vertices of $\Gamma$ are adjacent).
 This construction shows that we can work indifferently with
 1-vertex triangulations or with oriented maximal Wicks forms.

 Similarly, cellular decompositions of oriented surfaces with
 one vertex and one face correspond to oriented minimal 
 Wicks forms and were enumerated in \cite{[CM]}. The dual of an
 oriented minimal Wicks form is again a (generally non-equivalent)
 oriented minimal Wicks form and taking duals yields hence an
 involution on the set of oriented minimal Wicks forms.

 A vertex $V$ (with oriented edges $a,b,c$ pointing toward $V$) is
 {\em positive\/} if
 $$w=ab^{-1}\dots bc^{-1}\dots ca^{-1}\dots \quad {\rm or }\quad
 w=ac^{-1}\dots cb^{-1}\dots ba^{-1}\dots $$
 and $V$ is {\em negative\/} if    
 $$w=ab^{-1}\dots ca^{-1}\dots bc^{-1}\dots \quad {\rm or }\quad 
 w=ac^{-1}\dots ba^{-1}\dots ab^{-1}\dots \quad
 ..$$

 The {\em automorphism group\/} ${\rm Aut}(w)$ of an oriented
 Wicks form $w=w_1w_2\dots w_{2e}$ of length $2e$ is the group of all
 cyclic permutations $\mu$ of the linear word $w_1w_2\dots w_{2e}$ such
 that $w$ and $\mu(w)$ are isomorphic linear words (i.e. $\mu(w)$ is
 obtained from $w$ by permuting the letters of the alphabet). The group
 ${\rm Aut}(w)$ is a subgroup of the cyclic group ${\bf Z}/2e{\bf Z}$
 acting by cyclic permutations on linear words representing $w$.

 The automorphism group ${\rm Aut}(w)$ of an oriented Wicks
 form can of course also be described in terms of permutations on the
 oriented edge set induced by orientation-preserving
 homeomorphisms of $S$ leaving $\Gamma$ invariant. In particular
 an oriented maximal Wicks form and the associated dual
 1-vertex triangulation have isomorphic automorphism groups.

 We define the {\em mass\/} $m(W)$ of a finite set $W$ of oriented
 Wicks forms by
 $$m(W)=\sum_{w\in W}{1\over \vert{\rm Aut}(w)\vert}\quad .$$

 Let us introduce the sets
\par $W_1^g$: all oriented maximal
 Wicks forms of genus $g$ (up to equivalence),
\par $W^g_2(r)\subset W_1^g$: all oriented 
 maximal Wicks forms having an automorphism of order $2$ leaving 
 exactly $r$ edges of $w$ invariant by reversing their
 orientation. (This automorphism is the half-turn with respect to the
 \lq\lq midpoints" of these edges and exchanges the two adjacent
 vertices of an invariant edge.)
\par $W^g_3(s,t)\subset W_1^g$: all oriented maximal 
 Wicks forms having an automorphism of order $3$ leaving
 exactly $s$ positive and $t$ negative vertices invariant
 (this automorphism permutes cyclically the edges around
 an invariant vertex).
\par $W^g_6(3r;2s,2t)=W^g_2(3r)\cap W^g_3(2s,2t)$: all oriented
 maximal Wicks forms having an automorphism $\gamma$ of order 6
 with $\gamma^3$ leaving $3r$ edges invariant and $\gamma^2$
 leaving $2s$ positive and $2t$ negative vertices invariant
 (it is useless to consider the set $W_6^g(r';s',t')$ defined
 analogously since $3$ divides $r'$ and $2$ divides $s',t'$ if
 $W_6^g(r';s',t')\not=\emptyset$).

\par We define now the {\em masses\/} of these sets as
  $$\matrix
   {m_1^g\hfill&=&\displaystyle \sum_{w\in W_1^g}
    {1\over \vert{\rm Aut}(w)\vert}\quad ,\hfill \cr
    m_2^g(r)\hfill&=&\displaystyle \sum_{w\in W_2^g(r)}
    {1\over \vert{\rm Aut}(w)\vert}\quad ,\hfill \cr
    m_3^g(s,t)\hfill&=&\displaystyle \sum_{w\in W_3^g(s,t)}
    {1\over \vert{\rm Aut}(w)\vert}\quad ,\hfill \cr
    m_6^g(3r;2s,2t)\hfill&=&\displaystyle\sum_{w\in W_6^g(3r;2s,2t)}
    {1\over \vert{\rm Aut}(w)\vert}\quad .\hfill }$$

\par {\bf Theorem 1.1.}
{\sl\ \ (i) The group ${\rm Aut}(w)$ of automorphisms of an oriented
 maximal Wicks form $w$ is cyclic of order $1,\ 2,\ 3$ or $6$.
\par \ \ (ii) $\displaystyle
 m_1^g={2\over 1}\Big({1^2\over 12}\Big)^g{(6g-5)!\over g!(3g-3)!}
 \quad .$
\par \ \ (iii) $m_2^g(r)>0$ (with $r\in {\bf N}$) if and only if
 $f={2g+1-r\over 4}\in \{0,1,2,\dots\}$ and we have then
\par $\displaystyle m^g_2(r)={2\over 2}\Big({2^2\over 12}\Big)^f
 {1\over r!}{(6f+2r-5)!\over f!(3f+r-3)!}
 \quad .$
\par \ \ (iv) $m_3^g(s,t)>0$ if and only if
 $f={g+1-s-t\over 3}\in \{0,1,2,\dots\}$, $s\equiv 2g+1\pmod 3$ and
 $t\equiv 2g\pmod 3$ (which follows from the two previous conditions).
 We have then
\par $\displaystyle m^g_3(s,t)={2\over 3}\Big({3^2\over 12}\Big)^f
 {1\over s!t!}{(6f+2s+2t-5)!\over f!(3f+s+t-3)!}$ if $g>1$ and 
$\displaystyle m^1_3(0,2)={1\over 6}$.
\par \ \ (v) $m_6^g(3r;2s,2t)>0$ if and only if
 $f={2g+5-3r-4s-4t \over 12}\in \{0,1,2,\dots\}$,
 $2s\equiv 2g+1\pmod 3$ and $2t\equiv 2g\pmod 3$ (follows in fact
 from the previous conditions). We have then
\par $\displaystyle
 m_6^g(3r;2s,2t)={2\over 6}\Big({6^2\over 12}\Big)^f
 {1\over r!s!t!}{(6f+2r+2s+2t-5)!\over f!(3f+r+s+t-3)!}$
if $g>1$ and $\displaystyle m^1_6(3;0,2)={1\over 6}$.}

\par Set
 $$\matrix{m_2^g&=&\sum_{r\in {\bf N},\ (2g+1-r)/4\in {\bf N}\cup \{0\}}
 m_2^g(r)\quad ,\hfill \cr
 m_3^g&=&\sum_{s,t\in {\bf N},\ (g+1-s-t)/3\in {\bf N}\cup \{0\},\ 
 s\equiv 2g+1 \pmod 3} m_3^g(s,t)\quad ,\hfill \cr
 m_6^g&=&\sum_{r,s,t\in {\bf N},\ (2g+5-3r-4s-4t)/12\in
 {\bf N}\cup \{0\},\ 2s \equiv 2g+1\pmod 3} m_6^g(3r;2s,2t)\hfill}$$
 (all sums are finite) and denote by $M_d^g$ the number of equivalence
 classes of oriented maximal genus $g$ Wicks forms having an
  automorphism of order $d$ (i.e. an automorphism group with order
 divisible by $d$).

\par {\bf Theorem 1.2.}
 {\sl We have
 $$\matrix{ M_1^g=m_1^g+m_2^g+2m_3^g+2m_6^g\quad ,\hfill \cr
 M_2^g=2m_2^g+4m_6^g\quad ,\hfill \cr
 M_3^g=3m_3^g+3m_6^g\quad ,\hfill \cr
 M_6^g=6m_6^g \hfill} $$
 and $M_d^g=0$ if $d$ is not a divisor of $6$.}

\par The number $M_1^g$ of this Theorem is of course the number of
 inequivalent oriented maximal Wicks forms
 of genus $g$. The first 15 values $M_1^1,\dots,M_1^{15}$ are
 displayed in the Table at the end of this paper.

\par The following result is an immediate consequence of Theorem 1.2.

\par {\bf Corollary 1.1.}
 {\sl There are exactly\par $M_6^g$ inequivalent Wicks forms with $6$
 automorphisms,\par $M_3^g-M_6^g$ inequivalent Wicks forms with $3$
 automorphisms,\par $M_2^g-M_6^g$ inequivalent Wicks forms with $2$
 automorphisms and\par $M_1^g-M_2^g-M_3^g+M_6^g$ inequivalent
 Wicks forms without non-trivial automorphisms.}

\par {\bf Remark.}  Computing masses amounts to enumerating
 pointed objects, i.e. linear words instead of cyclic words
 in Definition 1.1. Their number is $(12g-6)m_d^g$, where
 $d$ is $ 1,2,3$ or $6$.

\par Let us remark that formula (ii) can be obtained from
 {\bf [WL]} (formula (9) on page 207 and the formula on the top of page
 211) or from {\bf [GS]} (Theorem 2.1 with $\lambda=2^{6g-3}$
 and $\mu=3^{4g-2}$).  We will reprove it independently.
 Related objects have also been considered in {\bf [HZ]}.


\section{Oriented Wicks forms}
\label{sec:2}

\par Let $V$ be a negative vertex of an oriented maximal Wicks
form of genus $g>1$. There are three possibilities, denoted
configurations of type $\alpha,\ \beta$ and $\gamma$
(see Figure 1) for the local configuration at $V$.

\centerline{\epsfysize2.5cm\epsfbox{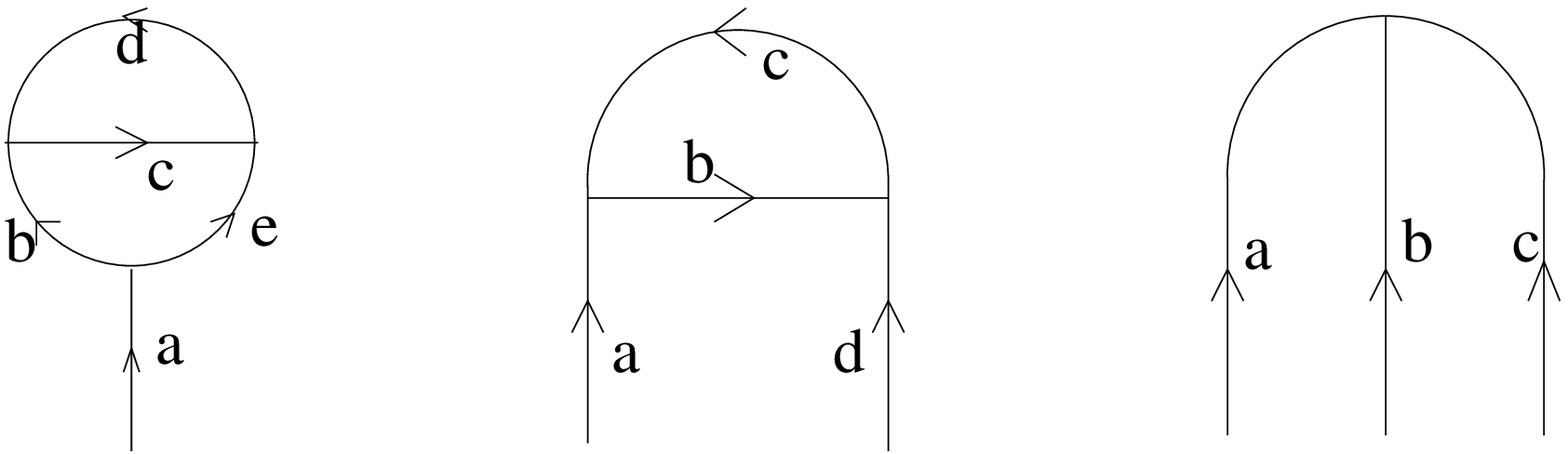}}
\centerline{\it Figure 1.}

\par Type $\alpha$. The vertex $V$ has only two neighbours 
which are adjacent to each other. This implies that $w$ is of the
form
$$w=x_1abcdb^{-1}ec^{-1}d^{-1}e^{-1}a^{-1}x_2u_1x_2^{-1}x_1^{-1}u_2$$
(where $u_1,u_2$ are subfactors of $w$)
and $w$ is obtained from the maximal oriented Wicks form 
$$w'=xu_1x^{-1}u_2$$ 
of genus $g-1$ by the substitution  $x\longmapsto
x_1abcdb^{-1}ec^{-1}d^{-1}e^{-1}a^{-1}x_2$ and $x^{-1}\longmapsto x_2^{-1}
x_1^{-1}$ (this construction is called the $\alpha -$construction in
{\bf [V]}).

\par Type $\beta$. The vertex $V$ has two non-adjacent neighbours. 
The word $w$ is then of the form
$$w=x_1abca^{-1}x_2u_1y_1db^{-1}c^{-1}d^{-1}y_2u_2$$
(where perhaps $x_2=y_1$ or $x_1=y_2$, see {\bf [V]} for all the details).
The word $w$ is then obtained by a $\beta -$construction from the 
word $w'=xu_1yu_2$ which is an oriented maximal Wicks form 
of genus $g-1$.

\par Type $\gamma$. The vertex $V$ has three distinct neighbours. We have
then
$$w=x_1ab^{-1}y_2u_1z_1ca^{-1}x_2u_2y_1bc^{-1}z_2u_3$$
(some identifications among $x_i,\ y_j$ and $z_k$ may occur,
see {\bf [V]} for all the details) and the word $w$ is obtained by a 
so-called $\gamma-$construction from the word $w'=x\tilde u_2y\tilde
u_1z\tilde u_3$.

\par {\bf Definition 2.1.}
We call the application which associates to an oriented maximal Wicks
form $w$ of genus $g$ with a chosen negative vertex $V$ the oriented
maximal Wicks form $w'$ of genus $g-1$ defined as above the {\it reduction}
of $w$ with respect to the negative vertex $V$. 

\par An inspection of figure $1$ shows that reductions with respect to
vertices of type $\alpha$ or $\beta$ are always paired since two doubly
adjacent vertices are negative, of the same type 
($\alpha$ or $\beta$) and yield the same reductions.

\par The above constructions of type $\alpha,\ \beta$ and $\gamma$
can be used for a recursive construction of all 
oriented maximal Wicks forms of genus $g>1$. 

\par {\bf Definition 2.2.} Consider an oriented maximal Wicks form
$w=w_1\dots w_{12g-6}$. To any
edge $x$ of $w$ we associate a transformation of $w$ called the
{\it IH-transformation on the edge $x$}.
Geometrically, an IH-transformation amounts to contracting the
edge $x$ of the graph $\Gamma\subset S$ representing the oriented
maximal Wicks form $w$. 
This creates a vertex of degree $4$ which can be split
in two different ways (preserving planarity of the graph on $S$)
into two adjacent vertices of degree $3$: 
The first way gives back the original Wicks form and the 
second way results in the IH-transformation. Graphically, an
IH-transformation amounts hence to replace a (deformed) letter I 
(a topological neighbourhood of the edge $x\in \Gamma\subset S$)
by a (deformed) letter H.

\par More formally, one considers the two subfactors $axb$ and $cx^{-1}d$ 
of the (cyclic) word $w$. Geometric considerations and Definition 1.1
show that $b\not=a^{-1},c\not=b^{-1},d\not=a^{-1},d\not=c^{-1}$ and
$(c,d)\not=(a^{-1},b^{-1})$. 

\par According to the remaining possibilities we consider now the following
transformation:

\par \ \ Type 1. $c\not=a^{-1}$ and $d\not=b^{-1}$. This implies that
$d^{-1}a^{-1}$ and $b^{-1}c^{-1}$ appear as subfactors in the cyclic word 
$w$. The IH-transformation on the edge $x$ is then defined by the
substitutions
$$\matrix{axb&\longmapsto &ab\cr
cx^{-1}d&\longmapsto &cd\cr
d^{-1}a^{-1}&\longmapsto&d^{-1}ya^{-1}\cr
b^{-1}c^{-1}&\longmapsto&b^{-1}y^{-1}c^{-1}}$$
in the cyclic word $w$.

\par \ \ Type 2a. Suppose $c^{-1}=a$. This implies that $b^{-1}axb$ and
$d^{-1}a^{-1}x^{-1}d$ are subfactors of the cyclic word $w$. Define 
the IH-transformation on the edge $x$ by
$$\matrix{b^{-1}axb&\longmapsto& b^{-1}yab\hfill\cr
d^{-1}a^{-1}x^{-1}d&\longmapsto&d^{-1}y^{-1}a^{-1}d\quad .}$$

\par \ \ Type 2b. Suppose $d^{-1}=b$. Then $axba^{-1}$ and $cx^{-1}b^{-1}
c^{-1}$ are subfactors of the cyclic word $w$ and we define the IH-
transformation on the edge $x$ by
$$\matrix{axba^{-1}&\longmapsto&abya^{-1}\hfill\cr
cx^{-1}b^{-1}c^{-1}&\longmapsto&cb^{-1}y^{-1}c^{-1}\quad .}$$

\par {\bf Lemma 2.1.} {\sl (i) 
IH-transformations preserve oriented maximal 
Wicks forms of genus $g$. 

\par (ii) Two oriented maximal Wicks forms related by a IH-transformation
of type 2 are equivalent.}

\par Proof. This results easily by considering the effect of
an IH-transformation on the graph $\Gamma\subset S$.\hfill QED

\par {\bf Proposition 2.1.} {\sl An oriented maximal Wicks 
form of genus $g$ has exactly $2(g-1)$ positive and $2g$ negative 
vertices.}

\par {\bf Lemma 2.2.} {\sl An $\alpha$ or a $\beta$ construction increases
the number of positive and negative vertices by 2.}

\par The proof is easy.

\par {\bf Lemma 2.3.} {\sl The number of positive or negative vertices is
constant under IH-transformations.}

\par Proof of Lemma 2.3. The Lemma holds for IH-transformations of type 2
by Lemma 2.1 (ii). Let hence $w,w'$ be two oriented maximal Wicks forms
related by an IH-transformation of type 1 with respect to the edge
$x$ of $w$ respectively $y$ of $w'$. This implies that $w$ contains the
four subfactors
$$axb\quad ,\quad cx^{-1}d\quad ,\quad d^{-1}a^{-1}\quad ,\quad b^{-1}
c^{-1}$$
and $w'$ contains the subfactors
$$ab\quad ,\quad cd\quad ,\quad d^{-1}ya^{-1}\quad ,\quad b^{-1}y^{-1}
c^{-1}$$
in the same cyclic order and they agree everywhere else. It is hence
enough to check the lemma for the six possible cyclic orders of the
above subfactors.

\par One case is
$$\matrix {
w\hfill
&=&axbu&\dots& cx^{-1}d&\dots& d^{-1}a^{-1}&\dots& b^{-1}c^{-1}&\dots\quad ,
\hfill \cr
w'\hfill
&=&abu&\dots& cd&\dots& d^{-1}ya^{-1}&\dots& b^{-1}y^{-1}c^{-1}&\dots\quad ,
\hfill\quad . }$$
In this case the two vertices of $w$ incident in $x$ and the two vertices
of $w'$ incident in $y$ have opposite signs. All other vertices are
not involved in the IH-transformation and keep their sign and the Lemma
holds hence in this case.

\par The five remaining cases are similar and left to the reader.
\hfill QED

\medskip
\par Proof of Proposition 2.1. The result is true in genus 1 by inspection
(the cyclic word $a_1a_2a_3a_1^{-1}a_2^{-1}a_3^{-1}$ is the unique 
oriented maximal Wicks form of genus 1 and has two negative vertices.) 

\par Consider now an oriented maximal Wicks form $w$ of genus $g+1$. 
Choose an oriented embedded loop $\lambda$
of minimal (combinatorial) length in $\Gamma$. 

\par First case.
If $\lambda$ is of length 2 there are two vertices related by a 
double edge in $\Gamma$. This implies that they are negative and of type
$\alpha$ or $\beta$. The assertion of Proposition 2.1 holds hence for 
$w$ by Lemma 2.2 and by induction on $g$.

\par Second case.
We suppose now that $\lambda$ is of length $\geq 3$. The oriented
loop $\lambda$ turns either left or right at each vertex.
If it turns on the same side at two consecutive vertices $V_i$ and $V_{i+1}$
the IH-transformation with respect to the edge joining $V_i$ and $V_{i+1}$
transforms $w$ into a form $w'$ containing a shorter loop. By Lemma
2.2, the oriented maximal Wicks forms $w$ and $w'$ have the same number 
of positive (respectively negative)  vertices.
\par 
If $\lambda$ does not contain two consecutive vertices $V_i$ and $V_{i+1}$ 
with the above 
property (ie. if $\lambda$ turns first left, then right, then left etc.)
choose any edge $\{V_i,V_{i+1}\}$ in $\lambda$ and make an 
IH-transformation with respect to
this edge. This produces a form $w'$ which contains a loop $\lambda'$ of
the same length as $\lambda$ but turning on the same side at the
two consecutive vertices $V_{i-1},V_i$ or $V_{i+1},V_{i+2}$. 
By induction on the length of $\lambda$ we can
hence relate $w$ by a sequence of IH-transformation to an oriented
maximal Wicks form $\tilde w$ of genus $g+1$ containing a loop of length 2 
for which the result holds by the first case. The Wicks forms $w$ and
$\tilde w$ have of course the same number of positive (respectively
negative) vertices by Lemma 2.2. \hfill QED

\section{Proof of Theorem 1.1}
\label{sec:3}

\par Proof of Theorem 1.1. 
We prove the corresponding assertions for oriented maximal Wicks forms.
The translation in terms of 1-vertex triangulations is immediate.

\par Let $w$ be an oriented maximal Wicks form with an automorphism
$\mu$ of order $d$. Let $p$ be a prime dividing $d$. The automorphism
$\mu'=\mu^{d/p}$ is hence of order $p$. If $p\not= 3$ then $\mu'$
acts without fixed vertices on $w$ and Proposition 2.1 shows that $p$
divides the integers $2(g-1)$ and $2g$ which implies $p=2$. 
The order $d$ of $\mu$ is hence 
of the form $d=2^a3^b$. Repeating the above argument with the prime 
power $p=4$ shows
that $a\leq 1$. 

\par All orbits of $\mu^{2^a}$ on the set
of positive (respectively negative)
vertices have either $3^b$ or $3^{b-1}$ elements and this leads 
to a contradiction if $b\geq 2$. 
This shows that $d$ divides $6$ and proves (i).

\par Proof of (ii).
An element of $W^{g+1}_1$ (which designs the set of
equivalence classes of oriented maximal Wicks forms with genus $g+1$)
can be obtained by applying
an $\alpha,\ \beta$ or $\gamma$ construction to an element in
$W_1^g$.

There are respectively $2{6g-3\choose
1},\ 4{6g-3\choose 2}+4{6g-3\choose 1}$ and 
$8{6g-3\choose 3}+8(6g-3)(6g-4)+8{6g-3\choose 1}$
possibilities for these constructions starting with a 
given element in $W_1^g$. On the other
hand, Proposition 2.1 shows that we can construct $2(g+1)$  oriented 
maximal Wicks forms in $W_1^g$ by applying reduction with respect to
a negative vertex to a given element in $W_1^{g+1}$.
The numbers of such
\lq\lq augmentations" and \lq\lq reductions" coincide after
weighting with the correct coefficients. These weights have to take care of
automorphisms and the fact that type $\alpha$ and $\beta$ constructions
give rise to $2$ negative vertices with the same \lq\lq inverse". 
A careful analysis shows that 
$$\Big(4{6g-3\choose 1}+8{6g-3\choose 2}+8{6g-3\choose 1}+
8{6g-3\choose 3}+16{6g-3\choose 2}+8{6g-3\choose 1}\Big)m_1^g=2(g+1)
m_1^{g+1}$$
which simplifies to
$$2(6g+1)(6g-1)(2g-1)m_1^g=(g+1)m_1^{g+1}$$
and proves (ii) by induction since
the function
$$g\longmapsto 2{(6g-5)!\over 12^g\ g!(3g-3)!}$$
satisfies the same recursion and we have equality for $g=1$
(since $m_1^g={1\over 6}=2{1!\over 12}{1!\over 1!\ 0!}$).

\par Proof of (iii). First case: $r<2g+1$ and hence $f={2g+1-r\over 4}\geq 1$.
Let $w$ be an oriented maximal Wicks form of genus $g$
with an automorphism $\mu$ of order 2 reversing the orientation of exactly
$r$ edges. There are ${6g-3-r\over 2}$ orbits of (unoriented) edges 
not invariant under $\mu$. Consider the graph obtained 
by removing all $\mu-$invariant edges from the
quotient graph $\Gamma/\mu$. After removing leaves and vertices of degree
$2$ we get an oriented maximal Wicks form $\tilde w$ with ${6g-3-r\over 2}-
r={3(2g-r-1)\over 2}$ edges and hence of genus $f={2g+1-r\over 4}\geq 1$
(recall that an oriented maximal Wicks form of genus $f$ has $3(2f-1)$
edges).

More precisely, let $w$ be represented by the word $w_1w_2\dots w_{12g-6}$.
The subword $w_1w_2\dots w_{6g-3}$ contains exactly one representant
of each orbit for the action of $\mu$ on oriented edges. Remove from
the word $w_1\dots w_{6g-3}$ all letters $w_k$ with 
$w_{k+6g-3}=w_k^{-1}$ ( they correspond to edges reversed by 
$\mu$). The resulting word $w'$ has length $6g-3-r$ and has the property
that if $w_k$ appears in $w'$ then either $w_k^{-1}$ or $w_{k+6g-3}^{-1}$
appears exactly once in $w'$ also. Replacing $w_{k+6g-3}^{-1}$ by
$w_k^{-1}$ we get a word which satisfies (i) of Definition 2.1. Removing
from this word (and of the resulting ones) all cyclic subfactors of the
form $w_kw_k^{-1}$ we get a word $w''$ satisfying also condition (ii). 
Cancel $w_i$ and its inverse (or $w_j$ and its inverse)
if $w_iw_j$ and $w_j^{-1}w_i^{-1}$ both occur as cyclic subfactors.
This produces ultimately an oriented maximal Wicks form $\tilde w$. 
A counting argument shows that it has genus $f$. (A good way to understand
what happens is to write the word $w$ along two concentric circles related
by radial segments indexed by invariant edges, see figure 2 below
illustrating the following example).

{\it Example.} Consider the oriented maximal Wicks form 
$$w=abcdea^{-1}fb^{-1}e^{-1}ghc^{-1}f^{-1}ig^{-1}d^{-1}h^{-1}i^{-1}$$ 
of genus two. The form $w$ admits an automorphism $\mu$ of order two 
(acting for instance by 
$\mu(abcdea^{-1}fb^{-1}e^{-1})=ghc^{-1}f^{-1}ig^{-1}d^{-1}h^{-1}i^{-1}$).
We have $r=1$, $g=2$ and $f={4+1-1\over 4}=1$ since $c$ is the unique 
edge invariant by $\mu$.
The subword $abcdea^{-1}fb^{-1}e^{-1}$ contains exactly one representant 
of each orbit for the action of $\mu$ on oriented edges.
Removing the unique edge $c$ reversed by $\mu$ and replacing $f$ by $d^{-1}$
(since $\mu(f)=d^{-1}$) we get $abdea^{-1}d^{-1}b^{-1}e^{-1}$. We cancel
$d$ and its inverse, since both $bd$ and $d^{-1}b^{-1}$ occur
as cyclic subfactors. The resulting word $\tilde w=abea^{-1}b^{-1}e^{-1}$ 
is an oriented Wicks form of genus 1.

An oriented maximal Wicks form $\tilde w$ obtained in this way
carries an extra structure defined as follows.
Write the word $w$ counterclockwise 
along two concentric circles related by $r$ radial segments in 
such a way that radial segments correspond to edges invariant under $\mu$
(see figure 2 below for a hopefully explanatory example).

\centerline{\epsfysize5.5cm\epsfbox{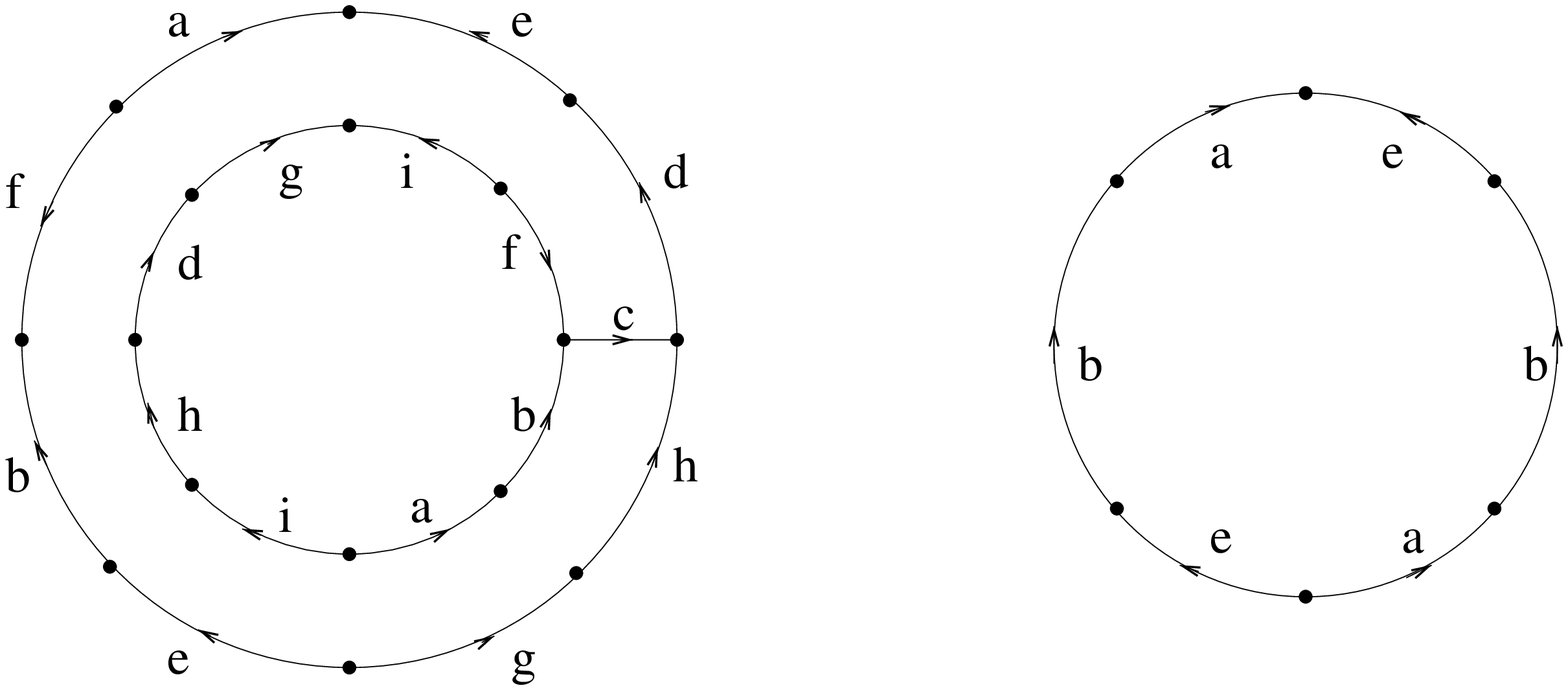}}
\centerline{\it Figure 2.}
\centerline{\it $w=abcdea^{-1}fb^{-1}e^{-1}ghc^{-1}f^{-1}ig^{-1}d^{-1}h^{-1}i^{-1}$}
\centerline{\it and the associated word $\tilde w=abea^{-1}b^{-1}e^{-1}$.}

\smallskip
Given a letter $\tilde l$ of $\tilde w$ choose a preimage $l$ of 
$\tilde l$. Set $\varphi(\tilde l)=0$ if $l$ and $l^{-1}$
are on the same circle and set $\varphi(\tilde l)=1$ otherwise
(we have hence $\varphi(a)=\varphi(b)=1$ and $\varphi(e)=0$
for the word $\tilde w=abea^{-1}b^{-1}e^{-1}$ of figure 2).
One checks then that $\varphi$ is well-defined and satisfies
$$\varphi(a)+\varphi(b)+\varphi(c)\equiv 0\pmod 2$$ 
whenever $a,b,c$ are $3$ edges incident in a common vertex of the
graph $\tilde G$ (which we identify of course with the word $\tilde w$).
Such a function is  called a ${\bf Z}/2{\bf Z}-${\it flow} 
on the graph $\tilde \Gamma$. 
\par Conversely, given an oriented maximal Wicks form $\tilde w$
of genus $f={2g+1-r\over 4}$ and a ${\bf Z}/
2{\bf Z}-$flow $\varphi$  on its graph $\tilde \Gamma$, 
we can construct       
$${(12f-6)(12f-2)\cdots(12f-10+4r)\over r!}$$
oriented maximal Wicks forms of genus $g$ having an automorphism $\mu$
of order 2 reversing $r$ edges associated to the pair $(\tilde w,\varphi)$. 
Indeed, we have 
$(12f-6)$ possibilities to attach the first edge reversed by $\mu$,
$(12f-2)$ choices for the second edge and so on. Since there are $r!$
possible orderings of the $\mu-$invariant edges we have to divide by $r!$.
Finally, the ${\bf Z}/2{\bf Z}$ flow shows
how to glue together preimages of orbits under $\mu$. 

\par The set of ${\bf Z}/2{\bf Z}-$flows is 
a vector space over ${\bf Z}/2{\bf Z}$ of dimension 
$2f$. This implies that we have
$$2^{2f}{(12f-6)(12f-2)\cdots(12f-10+4r)\over r!} m_1^f =2\ m_2^g(r) $$
(the factor $2$ on the right hand side comes from the fact that 
the Wicks forms contributing to $m_1^f$ are essentially weighted with
weight $1$ while they have weight $1\over 2$ in $m_2^g$).
This equation is also satisfied by replacing $m_1^f$ with $2{(6f-5)!\over
12^f f!(3f-3)!}$ and $m_2^g(r)$ with               
${(6f+2r-5)!\over 3^fr!f!(3f+r-3)!}$       
(recall that $g={4f+r-1\over 2}$) and this proofs (iii) in the first
case.

\par Second case: $f=0$ (the
construction of $\tilde w$ as above shows that we cannot have $f<0$). 
The idea is the same as in the first case.
Here we have to glue a first invariant edge on an empty word (1 
possibility) for the second and the third invariant edge we have 
2 possibilities, for the forth there are 6 possibilities etc.
Since there are no flows on an empty graph we get
$$2m_2^g(2g+1)=2{2\cdot 6\cdots (4r-10)\over r!}$$
which is readily checked.

\par Proof of (iv). First case: $t>0$. Let $w$ be an oriented maximal 
Wicks form having an automorphism of order $3$ fixing $s$ 
positive and $t>0$ negative vertices.
The $t$ fixed negative vertices give rise to $t$ possible reductions producing oriented Wicks forms $w'$ of genus
$g-1$ invariant under an automorphism of order $3$. The parameters 
of $w'$ are
then $(t-1,s)$. On the other hand, for any given oriented Wicks form $w'$
of genus
$g-1$ with an automorphism of order $3$ and parameters $(t-1,s)$ there
are $2(2g-3)$ $\gamma-$constructions yielding a Wicks form of genus $g$
with an automorphism of order $3$ and parameters $(s,t)$ (choose the
midpoint of any of the $6(2g-3)$ oriented edges in $w'$ and make the 
$\gamma-$construction with respect to its orbit).
We have hence
$$2(2g-3)m_3^{g-1}(t-1,s)=tm_3^g(s,t)$$
which is also satisfied by the righthand side of formula (iii) in Theorem
1.1.

\par Let us now consider the case $t=0$ (no invariant vertices of 
negative type). The proof of this case is very similar to the proof 
of (iii) .
\par We can suppose $g>1$ since  there are only two vertices 
 of negative type in genus 1 .
 We consider hence an oriented maximal Wicks form
$w$ of genus $g$ with an automorphism $\mu$ of order 3 fixing $s$ 
positive and no negative vertices. Since $\mu$ leaves no edge invariant,
there are ${6g-3\over 3}=2g-1$ orbits of edges. In genus $g>1$,
invariant vertices under an automorphism $\mu$ of order $3$ are
never adjacent. There are hence $s$ orbits of edges of $w$ incident in
a vertex fixed under $\mu$. Removing their orbits from the orbits of 
edges leaves us with a graph on the orbit space which has $s$ vertices of
degree $2$. Removing these vertices of degree $2$ yields an oriented
maximal Wicks form $\tilde w$ of genus $f={g+1-s\over 3}$ (all vertices
are of degree $3$, there is one face and there are $2g-1-2s=6f-3$ edges).
The construction of this form is completely analogous to the construction
in the proof of (ii). Graphically, one considers $3$ concentric
circles (indexed by the elements of ${\bf Z}/3{\bf Z}$) connected 
together by $s$ \lq\lq radial edges'' (which represent invariant
positive vertices together with their edges).
 
As in the proof of (ii) this form has an extra structure. 
This 
extra structure is here a ${\bf Z}/3{\bf Z}-$flow, ie. an application
$\varphi$ of the set of oriented edges of $\tilde w$ into 
${\bf Z}/3{\bf Z}$ such that $\varphi(e)\equiv -\varphi(-e)\pmod 3$ 
and $\varphi(a)+\varphi(b)+\varphi(c)\equiv 0\pmod 3$ 
for three oriented edges $a,b,c$ pointing
toward a common vertex of $\tilde w$.      

\par Conversely, given an oriented maximal Wicks form $\tilde w$ of
genus $f$ together with the above extra structure 
(a ${\bf Z}/3{\bf Z}-$flow on its graph $\tilde \Gamma$) there are
$$(12f-6)(12f-2)\cdots(12f-10+4s)\over s!$$
possibilities to \lq\lq extend" it into an oriented maximal Wicks form
$w$ of genus $g$ which has an automorphism $\mu$ of order 3 fixing
exactly $s$ positive and no negative vertices. 

\par Since the set of ${\bf Z}/3{\bf Z}-$flows 
on $\tilde \Gamma$ is a ${\bf Z}/3{\bf Z}-$vector
space of dimension $2f$ we get
$$3^{2f}{(12f-6)(12f-2)\cdots (12f-10+4s)\over s!}m_1^f=3\ m_3^g(s,0)
=3\ m_3^{3f+s-1}(s,0)\quad .$$
A routine calculation shows that this equation is also satisfied 
with $m_1^f$ replaced by $2{(6f-5)!\over 12^ff!(3f-3)!}$ and $m_3^{3f+s-1}
(s,0)$ replaced by
${2\over 3}\Big({3\over 4}\Big)^f{(6f+2s-5)!
\over s!(3f+s-1)!f!}$
and this proves (iv) in the case $f\geq 1$. The proof for $f=0$ is similar
to the analogous proof of (iii).

\par Proof of (v). We apply again the idea used in the proof of (iii). 
Let $w$ be an oriented maximal Wicks form with an automorphism $\mu$ of
order $6$. Considering the automorphism $\mu^3$ of order $2$ and
applying the reduction used in the proof of (iii) we get an oriented
maximal Wicks form $\tilde w$ of genus $h={2g+1-3r\over 4}$ together with
a ${\bf Z}/2{\bf Z}-$flow $\varphi$ on $\tilde \Gamma$. The Wicks form
$\tilde w$ is however an element of $W_3^h(s,t)$ and has hence
an automorphism $\tilde \mu$ of order $3$ which leaves $\varphi$
invariant. Analogously to the proof of (iii) 
we use this data to produce elements in $W_6^g(3r;2s,2t)$ by 
making all constructions $\tilde \mu-$invariant. We must understand
the vector space of $\tilde \mu-$invariant ${\bf Z}/2{\bf Z}-$flows:

\par {\bf Lemma 3.1.} 
{\sl Let $\tilde w\in W_3^h(s,t)$ be an oriented maximal 
Wicks form with an automorphism $\tilde \mu$ of order 3 (having parameters
$s,t$). The space of $\tilde\mu -$invariant ${\bf Z}/2{\bf Z}-$flows on
$\tilde \Gamma$ is then of dimension ${h+1-s-t\over 3}$.}

\par The lemma and a counting argument show then that
$$3^{(h+1-s-t)/3}{(4h-6)(4h-2)\cdots (4h-10+4r)\over r!}m_3^h(s,t)=
2m_6^g(3r;2s,2t)$$
and a routine calculation implies assertion (v).

\par Proof of Lemma 3.1. Let $\tilde \varphi $ be a $\tilde \mu-$invariant 
${\bf Z}/2{\bf Z}-$flow on $\tilde \Gamma$. We remark that $\tilde \varphi(a)
\equiv \tilde \varphi(b)\equiv \tilde \varphi(c)\equiv 0\pmod 2 $ if
$a,b,c$ are three edges incident in a $\tilde \mu-$fixed vertex. 
This shows that 
all reductions used in the proof of (iv) can also be applied to the flow 
$\tilde \varphi$ and these constructions are injective on 
$\tilde\mu-$invariant ${\bf Z}/2{\bf Z}-$flows. 

\par Theorem 1.1 is proved. \hfill QED

\par {\bf Table.}
The number of 1-vertex triangulations or of
oriented maximal Wicks forms in genus $1,\dots,15$:
$$
\matrix{
1\hfill &  1\hfill \cr 
2\hfill &  9\hfill \cr 
3\hfill &  1726\hfill \cr 
4\hfill &  1349005\hfill \cr 
5\hfill &  2169056374\hfill \cr 
6\hfill &  5849686966988\hfill \cr 
7\hfill &  23808202021448662\hfill \cr 
8\hfill &  136415042681045401661\hfill \cr 
9\hfill &  1047212810636411989605202\hfill \cr 
10\hfill &  10378926166167927379808819918\hfill \cr 
11\hfill &  129040245485216017874985276329588\hfill \cr 
12\hfill &  1966895941808403901421322270340417352\hfill \cr 
13\hfill &  36072568973390464496963227953956789552404\hfill \cr 
14\hfill &  783676560946907841153290887110277871996495020\hfill \cr 
15\hfill &  19903817294929565349602352185144632327980494486370\hfill \cr 
}$$


\begin{thebibliography}{99}

\bibitem{[B]} N.Biggs, {\em Algebraic Graph Theory}, Cambridge
 University Press, 1974 (Second Edition 1993).

\bibitem{[BC]} G.Bianchi, R.Cori, {\em Colorings of hypermaps
 and a conjecture of Brenner and Lyndon}, Pacific Journal of Maths,
 v.104, (1984), 41--048

\bibitem{[BI]} M. Bauer, C. Itzykson, {\em Triangulations\/} in {\em
 The Grothendieck Theory of Dessins d'Enfants\/} (Ed.L.Schneps), London
 Math. Soc. Lecture Notes Series 200, Cambridge University Press.

\bibitem{[BL]} J.L.Brenner, R.C.Lyndon, {\em Permutations
 and cubic graphs}, Pacific Journal of Maths, v.104, 285--315

\bibitem{[C]} M.Culler, {\em Using surfaces to solve equations
 in free groups}, Topology, v.20(2), 1981

\bibitem{[CCE]} J.A.Comerford, L.P.Comerford and C.C.Edmunds,
 {\em Powers as products of commutators},  Commun. in Algebra, 19(2),
 675--684 (1991)

\bibitem{[CE]} L.Comerford, C.Edmunds, {\em Products of commutators
 and products of squares in a free group}, Int.J. of Algebra
 and Comput., v.4(3), 469--480, 1994

\bibitem{[CM]} R.Cori, M.Marcus, {\em Counting non-isomorphic chord
 diagrams}, Theoretical Computer Science, 204, 55--73, 1998

\bibitem{[GS]} A.Goupil, G.Schaeffer, {\em Factoring
 n-cycles and counting maps of given genus}, Europ.J.Combinatorics,
 19(7), 819--834, 1998

\bibitem{[HZ]} J.Harer, D.Zagier, {\em The Euler Characteristic of
 the moduli space of curves}, Invent.Math, v.85, 457--486, 1986

\bibitem{[M]} L.Mosher, {\em A User's Guide to the Mapping Class
 Group: Once Punctured Surfaces}, DIMACS Series, Vol.25, 1994.

\bibitem{[V]} A.A.Vdovina, {\em Constructing Orientable Wicks Forms
 and Estimation of Their Number}, Communications in Algebra 23 (9), 
 3205--3222 (1995).

\bibitem{[WL]}T.R.S Walsh, A.B. Lehman, {\em Counting
 Rooted Maps by Genus}
 J.Combinatorial Theory, Ser.B, v.13, 192--218, 1972 

\end{thebibliography}
\end{document}